\documentclass[11pt]{article}

\usepackage[numbers, sort&compress]{natbib}
 \usepackage{setspace, tikz-cd}
\usepackage{amsmath,amssymb, amsthm, graphicx, wrapfig}
\usepackage[left=3cm,right=3cm, top=2.5cm,bottom=2.5cm,bindingoffset=0cm]{geometry}
\usepackage{multicol, mathtools}
\newtheorem{lemma}{Lemma}[section]
\newtheorem{statement}[lemma]{Proposition}
\newtheorem{remark-definition}[lemma]{Remark-Definition}
\newtheorem{theorem}[lemma]{Theorem}
\newtheorem{corollary}[lemma]{Corollary}

\newtheorem{proposition-conjecture}[lemma]{Proposition-conjecture}

\theoremstyle{definition}
\newtheorem{example}[lemma]{Example}

\newtheorem{definition}[lemma]{Definition}
\newtheorem{remark}[lemma]{Remark}

\mathtoolsset{showonlyrefs}                           

\newcommand{\Ker}[1]{\mathrm{Ker} \, #1}

\newcommand{\codim}[1]{\mathrm{codim} \, #1}
\newcommand{\corank}[1]{\mathrm{corank} \, #1}

\newcommand{\diff}[1]{\mathrm{d}  #1}

\newcommand{\diffFX}[2]{ \dfrac{\partial #1}{\partial #2} }

\newcommand{\diffX}[1]{ \frac{\partial }{\partial #1} }

\newcommand{\Z}{\mathbb{Z}}

\newcommand{\R}{\mathbb{R}}
\newcommand{\Complex}{\mathbb{C}}

\newcommand{\T}{\mathrm{T}}

\newcommand{\const}{\mathrm{const}}

\newcommand{\g}{\mathfrak{g}}

\newcommand{\GL}{\mathrm{GL}}

\newcommand{\calP}{\mathcal{P}}
\newcommand{\calQ}{\mathcal{Q}}

\title{Flat bi-Hamiltonian structures and invariant densities}
\author{Anton Izosimov\thanks{
Department of Mathematics,
University of Toronto, 
40 St. George Street, Toronto, Ontario, Canada M5S~2E4. 
e-mail: {\tt izosimov@math.utoronto.ca}}}
\date{}
\usepackage{caption, titlesec}

\titleformat{\section}{\large\bfseries\filcenter}{\thesection}{1em}{}
\titleformat{\subsection}{\bfseries}{\thesubsection}{1em}{}
\titleformat{\subsubsection}[runin]{\bfseries}{\thesubsubsection}{1em}{}[.]

\begin{document}
\maketitle
\begin{abstract}
A \textit{bi-Hamiltonian structure} is a pair of Poisson structures $\calP$, $\calQ$ which are \textit{compatible}, meaning that any linear combination $\alpha \calP + \beta \calQ$ is again a Poisson structure. A bi-Hamiltonian structure $(\calP, \calQ)$ is called \textit{flat} if $\calP$ and $\calQ$ can be simultaneously brought to a constant form in a neighborhood of a generic point.  We prove that a generic bi-Hamiltonian structure $(\calP, \calQ)$ on an odd-dimensional manifold is flat if and only if there exists a local density which is preserved by all vector fields Hamiltonian with respect to \nolinebreak $\calP$, as well as by all  vector fields Hamiltonian with respect to~$\calQ$.
\end{abstract}

%\tableofcontents
\section{Introduction}
Two Poisson structures $\calP$ and $\calQ$ on a manifold $M$ are called \textit{compatible} if any linear combination of them is again a Poisson structure. A pair of compatible Poisson structures is called a \textit{bi-Hamiltonian structure}. A vector field $v$ on $M$ is called \textit{bi-Hamiltonian} with respect to a bi-Hamiltonian structure $(\calP, \calQ)$ if it is Hamiltonian with respect to both $\calP$ and $\calQ$. \par
Bi-Hamiltonian structures play a fundamental role in the theory of integrable systems.
Since the pioneering works of F.\,Magri \nolinebreak \cite{Magri} and I.\,Gelfand and I.\,Dorfman \cite{GD}, it is known that bi-Hamiltonian systems automatically possess a large number of conservation laws and, as a consequence, tend to be completely integrable. Conversely, most of the known integrable systems possess a bi-Hamiltonian structure.  
 For this reason, bi-Hamiltonian structures have received a great deal of attention from the mathematical physics community.\par
In this note, we study local geometry of finite-dimensional bi-Hamiltonian structures. One of the most natural geometric questions about bi-Hamiltonian structures is whether there exists an analog of the Darboux theorem for compatible Poisson brackets. In other words, can two compatible Poisson tensors be simultaneously brought to a constant form in a neighborhood of a generic point? Besides being natural from the  point of view of abstract Poisson geometry, this problem also arises in the theory of integrable systems, in particular, in the context of separation of variables, see, e.g., \cite{FP2} and references therein.\par
\par
Geometry of bi-Hamiltonian structures has been intensively studied in the last 20 years. Main results in this field belong to I.\,Gelfand and I.\,Zakharevich \cite{GZ, GZ2, GZ3, Zakharevich}, P.J.\,Olver \cite{Olver},
 A.\,Panasyuk \nolinebreak \cite{Panasyuk}, and F.-J.\,Turiel \nolinebreak \cite{Turiel, Turiel2, Turiel4, Turiel3, Turiel5}.

 \par
 If two compatible Poisson tensors $\calP$ and $\calQ$ can be simultaneously brought to a constant form in a neighborhood of a generic point, then the bi-Hamiltonian structure $(\calP, \calQ)$ is called \textit{flat}. The main problem of the present paper is to find explicit conditions on tensors  $\calP$ and $\calQ$ under which the structure $(\calP, \calQ)$ is flat. We note that while flat bi-Hamiltonian structures are very non-generic objects, and one should expect flat examples to be rare, the situation is quite the opposite: most of bi-Hamiltonian structures naturally arising in mathematical physics are in fact flat. Here is just one remarkable example. 
 \begin{example}
 Let $\g$ be a finite-dimensional real Lie algebra, and let $a \in \g^*$. Then, on $\g^*$, there are two compatible Poisson structures, one being the standard \textit{Lie-Poisson bracket} $$\{f,g\}_\calP(x) := \langle x,[\diff f(x), \diff g(x)]\rangle,$$ and the other one being the \textit{frozen argument bracket} $$\{f,g\}_\calQ(x) := \langle a, [\diff f(x), \diff g(x)]\rangle.$$ For an arbitrary Lie algebra~$\g$, there is absolutely no reason for the bi-Hamiltonian structure~$(\calP, \calQ)$ to be flat. However, it is flat if $\g$ is \textit{simple} \cite{Panasyuk, Zakharevich}.

\end{example}

Note that, in bi-Hamiltonian geometry, one needs to distinguish between the even and odd dimensional cases. The reason for this comes from linear algebra. A generic skew-symmetric form on an even-dimensional vector space is non-degenerate. For two non-degenerate Poisson brackets \nolinebreak $\calP$ and \nolinebreak $\calQ$, one can define the \textit{recursion operator} $R = \calQ^{-1}\calP$. Clearly, if the structure $(\calP, \calQ)$ is flat, then the eigenvalues of $R$ must be constant. The converse is also true, as was proved by Turiel \cite{Turiel2}. \par

 In the odd-dimensional case, the situation is different. Any two generic pairs of forms on an odd-dimensional vector space are equivalent to each other. For this reason, generic odd-dimensional bi-Hamiltonian structures have no algebraic invariants, and the obstruction to flatness is of geometric nature.
 A fundamental theorem by Gelfand and Zakharevich \cite{GZ, GZ2} provides an isomorphism between the category of generic odd-dimensional bi-Hamiltonian structures and the category of \textit{Veronese webs}. As a corollary of the Gelfand-Zakharevich theorem, a bi-Hamiltonian structure is flat if and only if the associated Veronese web is \textit{trivializable}.\par
 The goal of the present note is to obtain a more intrinsic condition for flatness which does not appeal to the theory of webs. Our main observation is that the flatness problem for generic odd-dimensional bi-Hamiltonian structures is closely related to the notion of invariant densities, as introduced by A.\,Weinstein \cite{Weinstein}. Let $\calP$ be a Poisson structure on a manifold $M$. A density or, which is the same, a volume form on \nolinebreak $M$ is called \textit{invariant} if it is preserved by all vector fields Hamiltonian with respect to $\calP$. Likewise, if $M$ is endowed with a bi-Hamiltonian structure \nolinebreak $(\calP, \calQ)$, one can consider \textit{bi-invariant densities}, that are densities which are invariant with respect to both $\calP$ and $\calQ$. We say that a bi-Hamiltonian structure is \textit{unimodular} if it admits a bi-invariant density.\par
Our main result is that a generic odd-dimensional bi-Hamiltonian structure is flat if and only if it is locally {unimodular}. This result, albeit simple, is, to the best of our knowledge, first constructive criterion for flatness of bi-Hamiltonian structures. \par
 \medskip
\textbf{Acknowledgements.} The author is grateful to A.\,Bolsinov, A.\,Konyaev, I.\,Kozlov, and A.\,Oshemkov for fruitful discussions. 
This research was partially supported by the Dynasty Foundation Scholarship.
 \par\medskip
 \section{Pairs of bivectors on a vector space}
 In this section, we discuss some algebraic properties of pairs of bivectors on a finite-dimensional vector space. These properties serve as a motivation for the definition of a \textit{generic bi-Hamiltonian structure}, and also explain the difference between the even and odd dimensional cases in bi-Hamiltonian geometry.\par
Let $V$ be a finite-dimensional vector space. Consider the space $\Lambda^2V $ of bivectors on $V$ (i.e., bilinear skew-symmetric forms on $V^*$). Let
$
\Lambda_m := \{A \in \Lambda^2V \mid \corank A = m\}
$
be the set of bivectors of corank $m$. Then one has the following straightforward result.
\begin{statement}\label{codim}
The set $\Lambda_m$ is empty for $m \not\equiv \dim V \!\!\pmod 2$. For $m \equiv \dim V \!\!\pmod 2$, the set $\Lambda_m$ is a smooth submanifold of $\Lambda^2V$ of codimension 
$$
\codim \Lambda_m = \frac{1}{2} {m(m-1)}.
$$

\end{statement}

We say that a bivector is \textit{singular} if its rank is lower than maximum possible. From Proposition \nolinebreak \ref{codim} it follows that the set of singular bivectors has codimension $1$ if $\dim V$ is even, and codimension $3$ if $\dim V$ is odd. This fact is the reason for the difference between the even and odd dimensional cases in bi-Hamiltonian geometry.
\begin{definition}\label{linGenPair}
Let $V$ be an odd-dimensional vector space. We say that a pair of bivectors $A,B \in \Lambda^2V$ is \textit{generic} if $\alpha A + \beta B$ is non-singular for all $(\alpha, \beta) \in \Complex^2 \setminus \{0\}$, i.e. if
	$
 	\dim \Ker (\alpha A + \beta B) = 1 
 	$ for all $(\alpha, \beta) \in \Complex^2 \setminus \{0\}$. 
	\end{definition}
Let $\Omega = \Lambda^2V \times \Lambda^2V$ be the space of pairs of bivectors on $V$, and let $\Omega_0 \subset \Omega $ be the set of generic pairs. Then Proposition \ref{codim} implies the following result.

\begin{statement}\label{linGen}
 	$\Omega_0$ is an open dense subset of $\Omega$.
\end{statement}
\begin{proof}
The set of linear combinations $\alpha A + \beta B$ can be viewed as a projective line in the projectivization $ \mathbb P\Lambda^2V$ of the space of bivectors on $V$. The set $\Omega_0$ consists of pairs $(A,B) \in \Omega$ such that the line passing through $A$ and $B$ does not intersect (the projectivization of) the set of singular bivectors. Since the latter set has codimension $3$ in $ \mathbb P\Lambda^2V$, a generic line does not intersect it. The result follows.
\end{proof}
Furthermore, as follows from the proposition below, $\Omega_0$ is a homogeneous space for the natural $\GL(V)$-action on $\Omega$.  In other words, generic pairs of bivectors on an odd-dimensional vector space have no invariants. Note that in the even-dimensional case such invariants are given by the eigenvalues of the \textit{recursion operator} $A^{-1}B$.

\begin{statement}\label{canForm}
Let $V$ be a vector space of dimension $2n+1$. Assume that $(A,B) \in \Omega_0$. Then there exists a basis $e_1,  \dots, e_{2n+1}$ in $V$ such that
 	\begin{equation*}
 	A = \sum\nolimits_{i=1}^{n}e_i \wedge e_{i+n}, \quad
 	B = \sum\nolimits_{i=1}^{n} e_{i} \wedge e_{ i+n + 1}. 	
 	\end{equation*}
\end{statement}
\begin{proof}
The proof follows from the Jordan-Kronecker theorem; see, e.g., \cite{bolsinov2012jordan}.
\end{proof}
\par\medskip
\section{Generic bi-Hamiltonian structures on odd-dimensional manifolds}

In this section, we discuss the notions of \textit{generic} and \textit{flat} bi-Hamiltonian structures on an odd-dimensional manifold.

 \begin{definition}\label{genDef}
 	Let $M$ be an odd-dimensional manifold. A bi-Hamiltonian structure $(\calP, \calQ)$ on $M$ is called \textit{generic} at a point $z \in M$ if 
 	$
 	\dim \Ker (\alpha \calP(z) + \beta \calQ(z)) = 1$ for all  $(\alpha, \beta) \in \Complex^2 \setminus \{0\}$. 	
 \end{definition}
 In other words, $(\calP, \calQ)$ is {generic} at $z$ if the pair of bivectors $(\calP(z), \calQ(z))$ on $\T_{z}M$ is generic in the sense of Definition \ref{linGenPair}. From Proposition \ref{linGen} it follows that a bi-Hamiltonian structure generic at a point $z$ is also generic in a sufficiently small neighborhood of~$z$.

 \begin{remark}\label{Kron}
	In Definition \ref{genDef}, we follow the original terminology of Gelfand and Zakharevich \nolinebreak \cite{GZ,GZ2}. Note that generic bi-Hamiltonian structures on odd-dimensional manifolds are a particular case of \textit{Kronecker} structures. A bi-Hamiltonian structure $(\calP, \calQ)$ on a manifold $M$ is called Kronecker at a point $z \in M$ if the rank of $ \alpha \calP(z) + \beta \calQ(z)$ is the same for all $(\alpha, \beta) \in \Complex^2 \setminus \{0\}$.  
	So, a generic bi-Hamiltonian structure is the same as a Kronecker structure of corank one.\par
	
	 Kronecker bi-Hamiltonian structures are characterized by the following important property: if $(\calP, \calQ)$ is a Kronecker structure, then the set of all local Casimir functions of all brackets of the pencil $\alpha \calP + \beta\calQ$ is a completely integrable system; see \cite{Bolsinov}. In particular, this is so for generic bi-Hamiltonian structures on odd-dimensional manifolds.

\end{remark}

 \begin{definition}
 	A bi-Hamiltonian structure $(\calP, \calQ)$ is \textit{flat} in a neighborhood of a point $z$ if there exists a chart around $z$ in which both tensors $\calP$ and $\calQ$ have constant coefficients.
 \end{definition}
 By Proposition \ref{canForm}, any flat generic bi-Hamiltonian structure on $M^{2n+1}$ can be locally brought to the form
  	\begin{equation*}
 	\calP = \sum\nolimits_{i=1}^{n} \diffX{x_i} \wedge \diffX{x_{i+n}}, \quad
 	\calQ = \sum\nolimits_{i=1}^{n} \diffX{x_{i}} \wedge \diffX{x_{i+n + 1 }}.
 	\end{equation*}
\par\medskip
\section{Flatness via $\lambda$-Casimir families}
In this section, we recall the flatness criterion for generic odd-dimensional bi-Hamiltonian structures due to Gelfand and Zakharevich \cite{GZ3}.
    \begin{definition}
 	Let $(\calP, \calQ)$ be a bi-Hamiltonian structure on a manifold $M$. A family of functions  $F_\lambda \colon  M \to \R$ parametrized by $\lambda \in \R$ is called a $\lambda$-Casimir family if, for any $\lambda \in \R$, the function $F_\lambda$ is a Casimir function of $ \calP+\lambda \calQ$. A {$\lambda$-Casimir family} $F_\lambda$ is \textit{polynomial of degree $d$} if it is a degree $d$ polynomial in $\lambda$. A $\lambda$-Casimir family is \textit{non-singular} if $\diff F_\lambda \neq 0$ in $M$ for any $\lambda \in \R$.
 
 \end{definition}
  \begin{theorem}\label{thm0}
 	Let $(\calP, \calQ)$ be a bi-Hamiltonian structure on a manifold $M^{2n+1}$. Assume that $(\calP, \calQ)$ is generic at a point $z \in M^{2n+1}$. Then $(\calP, \calQ)$ is flat in a neighborhood of $z$ if 
	and only if	 near $z$ it admits a non-singular $\lambda$-Casimir family of degree $n$.
 \end{theorem}
 \begin{remark}
The condition of Theorem \ref{thm0} can be reformulated in terms of \textit{Veronese webs}; see~\cite{GZ3}. We do not use Veronese webs in the present paper.
\end{remark}
  \begin{remark}\label{degreeRem}
It can be shown that a generic bi-Hamiltonian structure in dimension $2n+1$ cannot have a non-singular $\lambda$-Casimir family of degree less than $n$. So, $n$ is actually the minimal possible degree of a $\lambda$-Casimir family. We also note that existence of a non-singular polynomial $\lambda$-Casimir family of degree higher than $n$ does not imply flatness. As an example, consider the following bi-Hamiltonian structure in $\R^3$.
\begin{align*}
\{x_1,x_2\}_{\calP} = 0, \quad \{x_3,x_1\}_{\calP} &= x_1, \quad \{x_3,x_2\}_{\calP} = -2x_2,\\
\{x_1,x_2\}_{\calQ} = 0, \quad \{x_3,x_1\}&_{\calQ} = 1, \quad \{x_3,x_2\}_{\calQ} = -2.
\end{align*}
This structure is generic at all points where $x_1 \neq x_2$ and has a $\lambda$-Casimir family of degree $3$ given by
$
F_\lambda= (x_1+\lambda)^2(x_2 + \lambda).
$
However, this structure is not flat; see Example \ref{nonFlat} below.
  \end{remark}

  Theorem \ref{thm0} was proved by Gelfand and Zakharevich in the analytic case \cite{GZ3}  and by Turiel in the smooth case  \cite{Turiel4}. Note that it is, in general, hard to use this result to prove or disprove flatness of an explicitly given bi-Hamiltonian structure. The problem is that any odd-dimensional bi-Hamiltonian structure admits infinitely many local $\lambda$-Casimir families, and it is not clear a priori whether there exists a distinguished one which is polynomial in \nolinebreak $\lambda$ and has degree $n$. Below, we show that $\lambda$-Casimir families of degree $n$ can be constructed by means of \textit{bi-invariant densities}. This allows us to give an alternative criterion for flatness which is easy to verify in particular examples; see Theorem~\ref{thm1} below.
\par\medskip
\section{Invariant densities on Poisson manifolds} 
In this section, we recall the notion of an \textit{invariant density}. For details, see \cite{Weinstein}.

\begin{definition} Let $M$ be a finite-dimensional manifold endowed with a Poisson structure $\calP$. A density (i.e., a non-vanishing top-degree form) $\omega$ on $M$ is \textit{invariant} with respect $\calP$ if is preserved by all vector fields which are Hamiltonian with respect to $\calP$. A Poisson structure which admits an invariant density is called \textit{unimodular}.
	\end{definition}
\begin{example}
Assume that $\calP$ is a non-degenerate Poisson structure on a manifold $M^{2n}$, i.e. $\calP^{-1}$ is a symplectic structure. Then the density
$
\omega := \Lambda^n(\calP^{-1})
$
is invariant with respect to $\calP$. Furthermore, any other $\calP$-invariant density coincides with $\omega$ up to a constant factor.
\end{example}
Note that any Poisson structure $\calP$ on $M$ is unimodular in a neighborhood of a generic point $z \in M$ (we say that a point $z$ is \textit{generic} for a Poisson structure $\calP$ if there exists a neighborhood of $z$ in which the rank of $\calP$ is constant). Indeed, by the Darboux theorem, there exists a coordinate system $(x_1, \dots, x_n)$  around $z$ in which the Poisson tensor $\calP$ has a constant form. In this coordinate system, an invariant density is given by
$
\omega := \diff x_1 \wedge \ldots \wedge \diff x_n
$.\par
 The proposition below gives a local description of all densities which are invariant with respect to a Poisson tensor $\calP$.
\begin{statement}\label{invCond}
Let $\calP$ be a Poisson structure on a manifold $M$. Let also $(x_1, \dots, x_n)$ be local coordinates on $M$, and let $\calP_{ij}$ be the components of the tensor $\calP$ in these coordinates. Then a density
$\omega = \rho(x_1, \dots, x_n)\diff x_1 \wedge \ldots \wedge \diff x_n$
is invariant with respect to $\calP$ if and only if
 	\begin{align}\label{invCondEq}
 		\sum\nolimits_{j=1}^n\left(\calP_{ij}\diffFX{\log \rho}{x_j} + \diffFX{\calP_{ij}}{x_j}\right) = 0.
 	\end{align}

\end{statement}

\begin{remark}\label{unimodularityInv}
The coordinate-free form of equation \eqref{invCondEq} is $\diff( \calP \star \omega) = 0$, where the star denotes contraction. In other words, an invariant density is the same as a top degree cycle in Poisson homology; see \cite{Weinstein}.
\end{remark}
\par\medskip

\section{Flatness via bi-invariant densities}
In this section, we formulate our main result: a necessary and sufficient condition for flatness of a generic odd-dimensional bi-Hamiltonian structure in terms of bi-invariant densities.\par
  \begin{definition}\label{biinv}
  A density $\omega$ is \textit{bi-invariant} with respect to a bi-Hamiltonian structure $(\calP, \calQ)$ if it is invariant with respect to both $\calP$ and $\calQ$. A bi-Hamiltonian structure $(\calP, \calQ)$ is called \textit{unimodular} if it admits a bi-invariant density.
 \end{definition}
 The following theorem is our main result.
 \begin{theorem}\label{thm1}
 	Let $(\calP, \calQ)$ be a bi-Hamiltonian structure on an odd-dimensional manifold $M$. Assume that $(\calP, \calQ)$ is generic at a point $z_0 \in M$. Then $(\calP, \calQ)$ is flat in a neighborhood of $z_0$ if 
	and only if	 it is unimodular  in a neighborhood of $z_0$.
 \end{theorem}

 The proof of Theorem \ref{thm1} is given in the next section. The following proposition shows how to apply this theorem in concrete examples.

  \begin{corollary}\label{cons}
 	Let $(\calP, \calQ)$ be a bi-Hamiltonian structure on $M^{2n+1}$ which is generic at a point $z_0 \in M^{2n+1}$. Let also $(x_1, \dots, x_{2n+1})$ be a coordinate system around $z_0$. Consider the following system of linear equations:
 	\begin{align}\label{eqns}
 	\begin{cases}
	\displaystyle
 		\,\vphantom{\int\limits_\int^\int}\sum\nolimits_{j=1}^{2n+1}\left(\calP_{ij}\alpha_j +\diffFX{\calP_{ij}}{x_j}\right) = 0,\\
 	\displaystyle	\vphantom{\int\limits_\int^\int}\,\sum\nolimits_{j=1}^{2n+1}\left(\calQ_{ij}\alpha_j +\diffFX{\calQ_{ij}}{x_j} \right)= 0,
 	\end{cases}
 	\end{align}
	where $\alpha_1, \dots, \alpha_{2n+1}$ are unknowns, and $\calP_{ij}, \calQ_{ij}$ are components of the Poisson tensors $\calP, \calQ$ in coordinates $(x_1, \dots, x_{2n+1})$.
 	Then the following statements hold.
 	\begin{enumerate}
 		\item If (\ref{eqns}) is solvable, then the solution is unique.
 		\item The bi-Hamiltonian structure $(\calP, \calQ)$ is locally unimodular if and only if (\ref{eqns}) is solvable, and the $1$-form \begin{equation}\label{alphaClosed}\alpha := \sum\nolimits_{i=1}^{2n+1} \alpha_i \diff x_i\end{equation} is closed.
 		\item If (\ref{eqns}) is solvable, and the form $\alpha$ given by~\eqref{alphaClosed} is closed, then
 		\begin{equation}\label{volumeForm}
 		\omega := \exp\left(\,\int\limits_{z_0}^{z}   \alpha \right) \diff x_1 \wedge \dots \wedge \diff x_{2n+1}
 		\end{equation}
 		is a bi-invariant density.
 	\end{enumerate}
 \end{corollary}
 \begin{proof}
 Statements 2 and 3 follow from Proposition \ref{invCond}. To prove the first statement, notice then the difference of any two solutions of (\ref{eqns}) belongs to $\Ker \cal P\, \cap \, \Ker \calQ$. At the same time, from Proposition \ref{canForm} it follows that $\Ker \calP \cap \Ker \calQ = 0$. \end{proof}
\begin{example}\label{nonFlat}
The bi-Hamiltonian structure considered in Remark \ref{degreeRem} is not flat.  Indeed, solving equations~\eqref{eqns}, we obtain
$$
 \sum\nolimits_{i=1}^{3} \alpha_i \diff x_i = \frac{\diff x_1}{x_2  - x_1} + \frac{\diff x_2}{2(x_2  - x_1)}.
$$
This form is not closed.
\end{example}
\begin{example}[Volterra lattice]
 Consider the bi-Hamiltonian structure of the periodic Volterra lattice; see, e.g.,~\cite{Damianou, FP}.
  This structure is given on $\R^{n}$ by
\begin{align*}
\begin{aligned}
 &\calP_{ij} := (\delta_{i+1,j} - \delta_{i,j+1})x_ix_j,\\
 \calQ_{ij} := (\delta_{i+1,j} - \delta_{i, j+1})&x_ix_j(x_i + x_j)+ \delta_{i+2,j}x_ix_{i+1}x_{i+2} -  \delta_{i-2,j}x_ix_{i-1}x_{i-2},
 \end{aligned}
 \end{align*}
 where all indices are modulo $n$ and $\delta_{ij}$ is the Kronecker delta. If $n$ is odd, this structure is generic almost everywhere. Let us prove that it is flat. The first of equations (\ref{eqns}) reads
$$
 		x_{i+1}\alpha_{i+1} - x_{i-1}\alpha_{i-1} = 0 \quad \forall\, i  \in \Z  /  n\Z.
$$
Since $n$ is odd, the latter equation implies 
$
x_i\alpha_i = x_j\alpha_j
$
for all $i,j \in  \Z  /  n\Z$. Denote $ x_i\alpha_i$ by~$\beta$. Then the second of equations (\ref{eqns}) takes the form
$$
\beta(x_i + x_{i+1}) - \beta(x_i + x_{i-1}) + \beta x_{i+1} - \beta x_{i-1} + 3(x_{i+1} - x_{i-1}) = 0,
$$
so
$\beta = -3 / 2$, and
$$
 \sum\nolimits_i \alpha_i \diff x_i = -\frac{3}{2}\sum\nolimits_i\frac{\diff{x_i}}{x_i}.
$$
This form is closed, so the structure ($\calP,\calQ$) is flat. The density
$$
\omega := \left(x_1 \ldots x_n\right)^{-3/2}\diff x_1 \wedge \dots \wedge \diff x_n
$$
is invariant with respect to $\calP$ and $\calQ$.
\end{example}
\begin{remark}
Let us comment on the dependence of the test proposed in Corollary \ref{cons} on the choice of coordinates. First, note that  if the system~\eqref{eqns} is consistent in one coordinate system, then it is also consistent in all coordinate systems. However, its solution $\alpha$ does depend on the choice of coordinates. 
A straightforward computation shows that solutions associated with coordinate systems $(x_i)$ and $( x'_i)$ are related by
\begin{align*}
\sum\nolimits_i \alpha_i  \diff x_i=  \sum\nolimits_{i} \alpha'_i\diff  x'_i + \diff \log \det J,
\end{align*}
where $J$ is the Jacobian of the transformation $(x_i) \to ( x'_i)$. From the latter formula it follows that the $2$-form
\begin{equation}\label{curvForm}
 \Theta := \diff \left( \sum \alpha_i  \diff x_i\right)
\end{equation}
does not depend on the choice of coordinates. In particular, if the $1$-form \eqref{alphaClosed} is closed in one coordinate system, then it is closed in all coordinate systems. Furthermore, if \eqref{alphaClosed} is closed,  then the density~\eqref{volumeForm} is well-defined up to a constant factor.\end{remark}

{
\begin{remark}\label{3dcase}
A straightforward computation shows that system~\eqref{eqns} is always consistent in dimension $3$. Therefore, in dimension $3$ the form \nolinebreak $\Theta$ given by~\eqref{curvForm} is always well-defined, and a generic bi-Hamiltonian structure on a $3$-manifold is flat if and only if $\Theta = 0$. It turns out that, up to a constant factor, the form $\Theta$ coincides with the curvature form defined in \cite{curv3}. This curvature form can be interpreted in the following two ways. Firstly, it coincides with the pullback of the Blaschke curvature of the $3$-web associated with the bi-Hamiltonian structure. Second, it is equal, up to a constant factor, to the skew-symmetric part of the Ricci tensor of any torsion-free affine connection $\nabla$ such that $\nabla \calP = \nabla \calQ = 0$. \par Apparently, there should exist a similar interpretation of the form $\Theta$ in any odd dimension. Note, however, that in higher dimensions system \eqref{eqns} is, in general, inconsistent. 

\end{remark}
\begin{remark}
Note that for Kronecker structures of higher corank (see Remark \ref{Kron} for the definition of Kronecker structures) unimodularity does not imply flatness (although it is still a necessary condition). As an example, consider two Poisson structures on $\R^4$ given by
\begin{align*}
\{x_1,x_2\}_{\calP} = x_2, \quad
\{x_1,x_3\}_{\calP}  &= x_3, \quad
\{x_1,x_4\}_{\calP}  = -2x_4,
\\
\{x_1,x_2\}_{\calQ} = 1, \quad
\{x_1,x_3\}_{\calQ}  &= 1, \quad
\{x_1,x_4\}_{\calQ}  = -2.
\end{align*}
It is easy to see that the bi-Hamiltonian structure $(\calP,\calQ)$ is Kronecker of corank $2$ at all points except those where $x_2 = x_3 = x_4$. Furthermore, $(\calP,\calQ)$ is a unimodular structure: a bi-invariant density  is given by
$
\diff x_1 \wedge \diff x_2 \wedge \diff x_3 \wedge \diff x_4.
$
 However, the structure $(\calP,\calQ)$ is not flat. To prove this, consider the intersection of kernels
 $
 \Ker \calP \cap \Ker \calQ.
 $
 The latter is a one-dimensional codistribution spanned by the $1$-form $$\beta = 2(x_4 - x_3)\diff x_2 + 2(x_2-x_4)\diff x_3 + (x_2-x_3)\diff x_4.$$
We have
 $
 \beta \wedge \diff \beta \neq 0,
 $
so the codistribution $ \Ker \calP \,\cap\, \Ker \calQ$ is not integrable. On the other hand, if the structure $(\calP,\calQ)$ was flat, this codistribution would be integrable. So, the bi-Hamiltonian structure $(\calP,\calQ)$ is not flat.
\end{remark}
\begin{remark}
After the first version of the present work became available on arXiv, Turiel published another criterion for flatness of generic bi-Hamiltonian structures:
\begin{theorem}[Turiel, \cite{turiel2015flatness}]
Let $(\calP, \calQ)$ be a bi-Hamiltonian structure on a manifold $M$ of dimension $2n+1 \geq 5$, and let $\Omega$ be an arbitrary density on $M$. Assume that $(\calP, \calQ)$ is generic at a point $z_0 \in M$. Then the bi-Hamiltonian structure $(\calP, \calQ)$ is flat near $z_0$ if and only if there exists a local $1$-form $\lambda$ such that \begin{align}\label{TC}\diff(\calP \star \Omega) = \lambda \wedge (\calP \star \Omega), \quad \diff(\calQ \star \Omega) = \lambda \wedge (\calQ \star \Omega).\end{align}
\end{theorem}
Let us show that condition \eqref{TC} is equivalent to unimodularity. First assume that the bi-Hamiltonian structure  $(\calP, \calQ)$ is unimodular, and let $\omega$ be a bi-invariant density. Let also $\Omega$ be an arbitrary density on $M$. Then, near $z_0$, we have $\Omega = f\omega$ for a certain non-vanishing function $f$, and
 $$
  \diff(\calP \star \Omega) = \diff (f (  \calP \star \omega)) =  \diff f \wedge (\calP \star \omega) = \diff \log f\wedge (\calP \star \Omega), %\quad   \diff(\calQ \star \Omega) = \diff f \wedge (\calQ \star \omega),
  $$
where in the second equality we used that $\omega$ is $\calP$-invariant, which, according to Remark~\ref{unimodularityInv}, is equivalent to $\diff(\calP \star \omega) = 0$. Analogously, we get $  \diff(\calQ \star \Omega) = \diff \log f\wedge (\calQ \star \Omega)$, so, taking $\lambda := \diff \log f$, we obtain \eqref{TC}.\par
Conversely, assume that  \eqref{TC} holds. Then, as shown by Turiel in the proof of Theorem 1 of \cite{turiel2015flatness}, the form $\lambda$ in~\eqref{TC} has to be closed, provided that $\dim M \geq 5$. But this means that locally one has $\lambda = \diff g$ for a certain smooth function $g$, and repeating the above computation we get that
$
\omega := \exp(-g)\Omega
$
is a bi-invariant density, as desired.
\end{remark}
\par\medskip
  \section{Proof of the flatness criterion}\label{3to2}
The proof of Theorem \ref{thm1} is based on the following lemma.
 \begin{lemma}\label{lemma2}
Let $\calP$ be a Poisson tensor of rank $2n$ on an odd-dimensional manifold $M^{2n+1}$. Let also $\omega$ be a density which is invariant with respect to $\calP$. Consider the $1$-form $\alpha$ obtained by contracting $\omega$ with the $n$th exterior power of $\calP$:
$
\alpha :=   \omega \star \Lambda^n\calP.
$
Then $\diff \alpha = 0$ and $\calP \alpha = 0$.
\end{lemma}
\begin{proof}
It suffices to prove the statement locally. By the Darboux theorem, there exists a local chart $(x_1, \dots, x_{2n+1})$ in which
$$
\calP = \sum_{i=1}^{n} \diffX{x_i} \wedge \diffX{x_{n+i}}.
$$
Let $\omega = \rho(x_1, \dots, x_n)\diff x_1 \wedge \ldots \wedge \diff x_n$. Then from Proposition \ref{invCond} and the invariance of $\omega$ it follows that $\calP \diff \rho = 0$, i.e. $\rho$ depends on $x_{2n+1}$ only:
$$
\omega = \rho(x_{2n+1})\,\diff x_1 \wedge \dots \wedge \diff x_{2n+1}.
$$
Therefore, we have 
$$
\alpha =  \omega \star \Lambda^n\calP =  \const \cdot \rho(x_{2n+1})\,\diff x_{2n+1},
$$
and the result follows.
\end{proof}

\begin{proof}[Proof of Theorem \ref{thm1}]
The implication ``flatness'' $\Rightarrow$ ``unimodularity'' is straightforward: if the tensors $\calP$ and $ \calQ$ are constant in a chart $(x_1, \dots, x_{2n+1})$, then the density
 $
\omega =  \diff x_1 \wedge \dots \wedge \diff x_{2n+1}
 $
 is invariant with respect to both $\calP$ and $\calQ$.\par 

Conversely, assume that there exists a density $\omega$ which is invariant with respect to $\calP$ and $\calQ$. Let
$$
\alpha_\lambda :=   \omega \star \Lambda^n( \calP +\lambda \calQ).
$$
This form is a degree $n$ polynomial in $\lambda$. By Lemma \ref{lemma2}, we have
$$
( \calP+\lambda \calQ)\alpha_\lambda = 0, \quad \diff \alpha_\lambda= 0.
$$
Also note that $\alpha_\lambda  \neq 0$ for any $\lambda \in \R$, since $\dim \Ker( \calP+\lambda \calQ) = 1$. Therefore, taking
$$
F_\lambda(z)  = \int_{z_0}^z \alpha_\lambda ,
$$
we obtain a non-singular polynomial $\lambda$-Casimir family of degree $n$. So, by Theorem \ref{thm0}, the bi-Hamiltonian structure $(\calP, \calQ)$ is flat, q.e.d.
\end{proof}

\bibliographystyle{plain}
\bibliography{flat}

\end{document}